\documentstyle[12pt]{article}
\renewcommand{\theequation}{\thesection.\arabic{equation}}
\newcounter{subequation}[equation]
\makeatletter

\expandafter\let\expandafter
\reset@font\csname reset@font\endcsname

\def\subeqnarray{\arraycolsep1pt
    \def\@eqnnum\stepcounter##1{\stepcounter{subequation}%
        {\reset@font\rm(\theequation\alph{subequation})}}
\jot5mm     \eqnarray}

\newfont{\blackb}{msbm10 scaled\magstep1}
\def\Bbb#1{\hbox{\blackb #1}}

\makeatother
\newcommand{\newsection}[1]{
\pagebreak[3]
\addtocounter{section}{1}
\setcounter{equation}{0}
\setcounter{subsection}{0}
\setcounter{footnote}{0}

\begin{flushleft}
{\Large\bf \thesection. #1}
\end{flushleft}
\nopagebreak
\medskip
\nopagebreak}
\newcommand{\newappendix}[1]{
\pagebreak[3]
\addtocounter{section}{1}
\setcounter{equation}{0}

\begin{flushleft}
{\Large\bf Appendix \thesection. #1}
\end{flushleft}
\nopagebreak
\medskip
\nopagebreak}

\begin{document}
\baselineskip 18pt
\begin{flushright}
FTUV 96-61 / IFIC 96-51
\end{flushright}
\vspace{1cm}
\begin{center}
\begin{Large}
{\bf Dirac operators on quantum $SU(2)$ group and quantum sphere}
\end{Large}
\\
\vspace{1cm}
{\bf P.N. Bibikov$^\dagger$ and P.P. Kulish$^{\dagger\ddagger}$}
\vskip10pt
{\it $^\dagger$ St.Petersburg Branch of Steklov Mathematical 
Institute}\vskip0pt
{\it Fontanka 27, St.Petersburg, 191011, Russia}
\\
\vskip10pt
{\it $^\ddagger$ Departamento de F\'{\i}sica Te\'orica and IFIC,}\vskip0pt
{\it Centro Mixto Universidad de Valencia-CSIC}\vskip0pt
{\it E-46100-Burjassot (Valencia) Spain.}
\end{center}

\begin{abstract}                
Definition of Dirac operators on the quantum group $SU_{q}(2)$
and the quantum sphere $S^{2}_{q \mu}$ are discussed. In both cases 
similar $SU_{q}(2)$-invariant form is obtained. It is 
connected with corresponding Laplace operators.
\end{abstract}
\setcounter{section}{-1}

\newsection{Introduction}

Dirac operator undoubtedly is one of the basic notions of the A. Connes
approach to noncommutative differential geometry \cite{1}. Therefore
it is natural to define Dirac operators on quantum groups and
quantum homogeneous spaces
as one of the most important and studied examples of quantum manifolds.
As it will be shown in this paper these two problems are closely connected.
For instance the same formula defines Dirac operator both on the quantum
group $SU_{q}(2)$ and on the quantum 2-sphere $S^{2}_{q \mu}$ \cite{3}.

In the paper \cite{4} the Dirac operator on the special case of quantum
sphere was proposed. The leading idea of this approach was to embed quantum
sphere inside quantum three dimensional Euclidean space and then to
construct in this space operator commuting with the radius of quantum sphere
and in the commutative limit coinciding with the standard Dirac operator on
the usual two dimensional sphere.

In the proposed approach differential structures on $SU_{q}(2)$ and
$S^{2}_{q\mu}$ are introduced according to the right $SU_{q}(2)$-coaction.
That is why Dirac operators on $SU_{q}(2)$ and $S^{2}_{q\mu}$ have
the similar form. The connection between \cite{4} and the proposed
approach will be studied in the next paper.

The $q$-Dirac operators appeared also in the study of the $q$-deformed
Poincare groups and $q$-Minkowski spaces (see a review \cite{13}). Our
covariance approach is similar to the one used in \cite{13}.

Recently the notions of noncommutative geometry were used to construct
noncommutative manifold started from standard sphere and the representation
theory of the $SU(2)$ group \cite{14}. The regularization parameter is
connected with the highest spin. It has to be pointed out that due to the
complete analogy of the ${\rm su}(2)$ algebra representation theory and the
quantum ${\rm su}_{q}(2)$ one and above mentioned $q$-Dirac operator,
the corresponding constructions can be extended to the quantum sphere case
giving rise to extra parameters in the theory.

The paper is organized as follows. In Sect. 2 we give a detailed notion
on the quantum group $SU_{q}(2)$ and the quantum sphere $S^{2}_{q \mu}$. 
Sect. 3 is
devoted to the notion of the quantum enveloping algebra ${\rm su}_{q}
(2)$. In Sect. 4 we give the notion of the algebra of functions on the
$SU_{q}(2)$ quantum cotangent bundle and the analogous definition of
the quantum 2-sphere cotangent bundle. And finally in Sect. 5
the $SU_{q}(2)$-covariant construction of Dirac operator is proposed.

\newsection{Quantum $SU_q(2)$ group and quantum 2-sphere $S^{2}_{q \mu}$ }
The algebra of functions on the quantum group 
$F_{q}(G) = {\rm Fun}(SU_q(2))$ is  an associative $^*$-algebra
generated by two elements $a$  and $b$ satisfying relations
\cite{2} ($\lambda \equiv q - \frac{\displaystyle 1}{\displaystyle q}$)
\begin{eqnarray}
ab  =  qba, \qquad   ab^*  =  qb^*a,  \qquad b^*b  = bb^* ,
\nonumber\\
a^*a - aa^*  =  \frac{1}{q}\lambda b^*b, \qquad  aa^* + b^*b = I \,.
\end{eqnarray}
These relations can be written in a compact matrix form \cite{2}
\begin{eqnarray}
\hat RT_{1} T_{2} & = & T_{1}T_{2} \hat R, \nonumber\\
{\rm det}_{q}T & = & aa^* + b^*b = I \,, 
\end{eqnarray}
where $ T_{1} = T \otimes I_{2}$, $T_{2} = I_{2} \otimes T $ and
$  I_{2}$  is the $ 2 \times 2 $ unit matrix,
\begin{equation}
T =
\left( \begin{array}{cc}
a & b \\
-\frac{\displaystyle 1}{\displaystyle q}b^* & a^*
\end{array} \right)
\end{equation}
is the standard matrix of generators of $F_{q}(G)$ and the $R$-matrix
\begin{eqnarray}
\hat R =
\left( \begin{array}{cccc}
q & & & \\
& \lambda & 1 & \\
& 1 & & \\
& & & q
\end{array} \right)
\end{eqnarray}
satisfies the Yang-Baxter equation (in the braid group form) 
\begin{equation}
(\hat R \otimes I_{2})(I_{2} \otimes \hat R)(\hat R \otimes
I_{2})  =  (I_{2}  \otimes \hat R)(\hat R \otimes
I_{2})(I_{2} \otimes \hat R)
\end{equation}
and the Hecke condition
\begin{equation}
{\hat R}^2 =   \lambda {\hat R} + I_4 \,,
\end{equation}
where $I_4$ is the $4 \times 4$ unit matrix.

Comultiplication
\begin{equation}
\triangle T = T ( \otimes ) T
\end{equation}
(i.e. $\triangle T_{ik} = \sum_{j} T_{ij} \otimes T_{jk} $ ), antipode
\begin{eqnarray}
S(T) =
\left( \begin{array}{cc}
a^* & -\frac{\displaystyle 1}{\displaystyle q}b \\
b^* & a
\end{array} \right)
\end{eqnarray}
and counit $\varepsilon (T) = I_2 $ define the structure of
Hopf algebra on $ F_{q}(G) $
\cite{2}. We have also a useful
relation 
\begin{equation} 
S(T)T = TS(T) = I_2 \,.
\end{equation}

The algebra of functions on quantum 2-sphere $F_{q}(S)={\rm Fun}(S^2_{q\mu})$
is an associative $^*$-algebra with three generators 
$ x_{+} \,, x_{-}\,, x_{3} $ and two parameters $q, \mu$\,,
satisfying the following relations \cite{3}
\begin{eqnarray}
x_{+}^{*} = x_{-},\qquad  x_{3}^{*} & = & x_{3}, \\
qx_{3}x_{+} - \frac{1}{q}x_{+}x_{3} & = & \mu{x_{+}}, \nonumber \\
\lambda x_{3}^{2} + \frac{1}{[2]_{q}}(x_{+}x_{-} - x_{-}x_{+})
& = & \mu{x_{3}}, \nonumber \\
qx_{-}x_{3} - \frac{1}{q}x_{3}x_{-} & = & \mu{x_{-}}, \nonumber \\
x_{3}^{2}+\frac{1}{[2]_{q}}(qx_{-}x_{+}+\frac{1}{q}x_{+}x_{-}) & = & r^{2}
\end{eqnarray}
where $r^{2}$ is the central element of $F_{q}(S)$. These relations can be 
written also in a compact matrix form \cite{8} 
\begin{eqnarray}
M = 
\left( \begin{array}{cc}
\frac{\displaystyle 1}{\displaystyle q}x_3 & x_ - \\
x_{+} & - qx_{3}
\end{array} \right)
\end{eqnarray}

\begin{eqnarray}
M^{\dagger} &=& M,   \\
{[}\hat R,(M_{2}\hat RM_{2} + \mu{q}M_{2}){]} & = & 0,
\nonumber  \\
\frac{1}{[2]_q}{\rm tr}_qM^2 & = & r^2
\end{eqnarray}
where $ [N]_q = \frac{\displaystyle 1}{\displaystyle \lambda}(q^N - q^{-N})$
and the $q$-trace for
arbitrary $ 2\times 2$ matrix $X$ is defined by  ${\rm tr}_{q}X =
{\rm tr}DX$ for $ D = {\rm diag}(q,\frac{\displaystyle 1}
{\displaystyle q})$. It has an invariance property
\begin{equation}
{\rm tr}_{q}S(T)XT = {\rm tr}_{q}X
\end{equation}
where $X$ is an arbitrary $2\times 2$ matrix, whose elements are commuting
with elements of $T$.

The relations (1.11), (1.14) can be represented in the form of 
the reflection equation \cite{7}, \cite{8}:
\begin{eqnarray}
\hat R{\Bbb M}_{2}\hat R{\Bbb M}_{2}& =& {\Bbb M}_{2}\hat R{\Bbb M}_{2}
\hat R, \nonumber \\
\frac{1}{[2]q}{\rm tr}_{q}{\Bbb M}^{2}& =& (\mu^{2}q^{2}+\lambda^{2}r^{2})
\end{eqnarray}
where ${\Bbb M} = \mu{q}I_{2} + \lambda M$, thus omitting the condition 
${\rm tr}_{q}M = 0 $.

The following right coaction of the quantum group 
$\varphi_{R}: F_{q}(S)\rightarrow F_{q}(S)
\otimes F_{q}(G)$
\begin{eqnarray}
\varphi_{R}({\Bbb M})& =& S(T){\Bbb M}T, \nonumber \\
(\varphi_{R}({\Bbb M}_{ij})& =& \sum {\Bbb M}_{kl}\otimes{S(T)_{ik}}T_{lj})
\end{eqnarray}
according to eqs. (1.2), (1.16) defines on $F_{q}(S)$ structure of 
$F_{q}(G)$-comodule algebra.

In the commutative case $q = 1$, $\mu = 0$ eqs. (1.12) give
\begin{equation}
x_{1}^{2} + x_{2}^{2} + x_{3}^{2} = r^{2}
\end{equation}
for $x_{1} = \frac{\displaystyle 1}{\displaystyle 2}(x_{+} + x_{-})$ and
$x_{2} = \frac{\displaystyle 1}{\displaystyle 2i}(x_{+} -x_{-})\,,$ 
while for $ \mu \neq 0 $ one gets the Lie algebra $sl(2)$. 

Invariant scalar products $\langle\cdot,\cdot\rangle_{G}$ on $F_{q}(G)$ and
$\langle\cdot,\cdot\rangle_{S}$ on $F_{q}(S)$ can be defined as follows
\cite{11},\cite{15}:
\begin{eqnarray}
\langle u_{G},v_{G}\rangle_{G} &=& h_{G}(u^*_{G}v_{G}), \nonumber \\
\langle u_{S},v_{S}\rangle_{S} &=& h_{S}(u^*_{S}v_{S}), \nonumber \\
u_{G},v_{G} \in F_{q}(G),&& u_{S},v_{S} \in F_{q}(S).
\end{eqnarray}
The map $h_{G}$ is the Haar measure on $F_{q}(G)$ \cite{11}, i.e. the positive 
linear functional $h_{G}: F_{q}(G) \rightarrow {\Bbb C}$ invariant under
the coproduct
\begin{eqnarray}
(id\otimes h_{G})\triangle(u_{G}) & = & h_{G}(u_{G})I,  \nonumber \\
(h_{G}\otimes id)\triangle(u_{G}) & = & h_{G}(u_{G})I,  \nonumber \\
u_{G}\in F_{q}(G). & &
\end{eqnarray}
and $h_{S}$ is the $\varphi_{R}$-invariant positive linear functional 
on $F_{q}(S)$ \cite{15}
\begin{eqnarray}
(h_{S}\otimes id)\varphi_{R}(u_{S}) & = & h_{S}(u_{S})I, \nonumber \\
u_{S}\in F_{q}(S)\,.
\end{eqnarray}

\newsection{Quantum universal enveloping algebra ${\rm su}_{q}(2)$}
Quantum universal enveloping algebra ${\rm su}_{q}(2)$ is 
an associative algebra
generated by four elements $k, k^{-1}, e, f$ and relations \cite{2},\cite{9}
\begin{eqnarray}
kk^{-1} & = & I,\qquad  k^{-1}k = I, \nonumber \\
ek & = & qke,\qquad kf = qfk, \nonumber \\
k^2 - k^{-2} & = & \lambda(fe - ef)
\end{eqnarray}
or in the matrix form \cite{2}
\begin{eqnarray}
\hat RL^{\pm}_{2}L^{\pm}_{1} & = & L^{\pm}_{2}L^{\pm}_{1}\hat R \,, 
\nonumber \\
\hat RL^{+}_{2}L^{-}_{1} & = & L^{-}_{2}L^{+}_{1}\hat R \,, \nonumber \\
{\rm det}_{q}L^{\pm} & = & I \,,
\end{eqnarray}
where
\begin{eqnarray}
L^+ =
\left(\begin{array}{cc}
k^{-1}&\frac{\displaystyle \lambda}{\displaystyle \sqrt{q}}f \\
0 & k
\end{array}\right), \, \, \,
L^{-} =
\left(\begin{array}{cc}
k&0 \\
-\lambda \sqrt{q}e&k^{-1}
\end{array}
\right)
\end{eqnarray}
and $\hat R$ is given by (1.4).

Comultiplication
\begin{eqnarray}
\triangle k = k\otimes k, & & \triangle k^{-1} = k^{-1}\otimes k^{-1},
\nonumber \\
\triangle e = e\otimes k + k^{-1}\otimes e, & & \triangle f =
f\otimes k + k^{-1}\otimes f
\end{eqnarray}
and counit map 
\begin{eqnarray}
\varepsilon(k) = 1, & & \varepsilon(k^{-1}) = 1, \nonumber \\
\varepsilon(e) = 0, & & \varepsilon(f) = 0
\end{eqnarray}
or in the compact matrix form $\triangle L^{\pm} = L^{\pm} 
(\otimes) L^{\pm}$ and
$\varepsilon(L^{\pm}) = I_2$
and antipode
\begin{eqnarray}
S(L^+) =
\left(\begin{array}{cc}
k&-\lambda\sqrt{q}f \\
0&k^{-1}
\end{array}\right), \qquad
S(L^{-})=
\left(\begin{array}{cc}
k^{-1}&0 \\
\frac{\displaystyle \lambda}{\displaystyle \sqrt{q}}e&k
\end{array}\right)
\end{eqnarray}
define the structure of Hopf algebra on  ${\rm su}_{q}(2)$. 
As in the $F_{q}(G)$ case we have 
relations
\begin{equation}
L^{\pm}S(L^{\pm}) = S(L^{\pm})L^{\pm} = I_{2}
\end{equation}

The Hopf algebres $F_{q}(G)$ and ${\rm su}_{q}(2)$ are related 
by the duality \cite{2}.
This means that there exists a pairing $\langle\cdot,\cdot\rangle: {\rm su}_{q}
(2)\otimes F_{q}(G) \rightarrow {\Bbb C}$ satisfying relations:
\begin{subeqnarray}
\langle t_{1}t_{2},v\rangle & = & \langle t_{1}\otimes t_{2},\triangle v
\rangle.  \\
\langle t,v_{1}v_{2}\rangle & = & \langle\triangle (t), v_{1}\otimes v_{2}
\rangle.  \\
\langle t,S(v)\rangle & = & \langle S(t),v\rangle.  \\
\langle I,v\rangle & = & \varepsilon(v).  \\
\langle t,I\rangle & = & \varepsilon(t).
\end{subeqnarray}
\begin{displaymath}
t,t_{1},t_{2}\in U, \qquad v,v_{1},v_{2}\in F_{q}(G).
\end{displaymath}
Using the matrix $T$ of generators the pairing between 
$F_{q}(G)$ and ${\rm su}_{q}(2)$ can be
defined by \cite{2}:
\begin{eqnarray}
\langle k,T\rangle = \left( \begin{array}{cc}
                  \frac{\displaystyle 1}{\displaystyle \sqrt{q}}&0 \\
                  0 & \sqrt{q}
               \end{array} \right), & &
\langle k^{-1},T\rangle = \left( \begin{array}{cc}
                       \sqrt{q} & 0 \\
                       0 & \frac{\displaystyle 1}{\displaystyle \sqrt{q}}
                    \end{array} \right), \nonumber \\
\langle e,T\rangle = \left( \begin{array}{cc}
                  0 & 1 \\
                  0 & 0
               \end{array} \right), & &
\langle f,T\rangle = \left( \begin{array}{cc}
                  0 & 0 \\
                  1 & 0
               \end{array} \right)
\end{eqnarray}
or in the compact matrix form \cite{2}
\begin{equation}
\langle L^{\pm}_{1},T_{2}\rangle = R^{\pm}
\end{equation}
where
\begin{equation}
R^{+}=\frac{1}{\sqrt{q}}\hat{R}{\cal P} ,\qquad  R^{-}=
\sqrt{q}\hat{R}^{-1}{\cal P}
\end{equation}
and
\begin{displaymath}
{\cal P} =
\left(\begin{array}{cccc}
1&&& \\
&&1& \\
&1&& \\
&&&1
\end{array}\right)
\end{displaymath}
is a $4\times{4}$ permutation matrix in ${\Bbb C}^2 \otimes {\Bbb C}^2$ , 
so that for every $2\times2$ matrix
$X: X_{2} ={\cal P}X_{1}{\cal P}$.

The center of ${\rm su}_{q}(2)$ is generated by the element
\begin{equation}
C_{q} = \frac{1}{q}k^{2} + qk^{-2} + \lambda^{2}fe \,.
\end{equation}

Inside the algebra ${\rm su}_{q}(2)$ there exists an important
subalgebra generated by elements of the matrix
\begin{equation}
{\Bbb L}=L^{+}S(L^{-}) \,.
\end{equation}
The matrix ${\Bbb L}$ satisfies the reflection equation 
( see e.g. \cite{5},\cite{7} )
\begin{equation}
\hat {R}{\Bbb L}_{2}\hat {R}{\Bbb L}_{2}={\Bbb L}_{2}\hat{R}{\Bbb L}_{2}
\hat{R} \,.
\end{equation}
The reflection equation (2.14) is invariant under 
the right $SU_{q}(2)$-coaction
\begin{equation}
\varphi_{R}({\Bbb L}) = S(T){\Bbb L}T \,,
\end{equation}
hence the entries of ${\Bbb L}$ are generators of a right 
$SU_{q}(2)$-comodule algebra. In terms of (2.1), (2.3) it is 
\begin{eqnarray}
{\Bbb L}= \left( \begin{array}{cc}
k^{-2} + \frac{\displaystyle 1}{\displaystyle q} \lambda^{2} f e &
\frac{\displaystyle 1}{\displaystyle \sqrt{q}} \lambda f k \\
\frac{\displaystyle 1}{\displaystyle \sqrt{q}} \lambda k e & k^{2}
\end{array} \right)
\end{eqnarray}
and the central element is 
\begin{equation}
C_{q}={\rm tr}_{q}{\Bbb L} \,. 
\end{equation}
According to (2.15), (2.17) and (1.15) this $C_{q}$ is the invariant
element of this coaction.
\begin{equation}
\varphi_{R} (C_{q}) = C_{q} \otimes {I} \,. 
\end{equation}

We can also represent matrix ${\Bbb L}$ in the form
\begin{equation}
{\Bbb L}=\frac{1}{[2]_{q}}C_{q}I_{2}+\frac{\lambda}{q}L_{q}
\end{equation}
with the traceless matrix $L_{q}$ , $tr_q L_{q} = 0$ 
\begin{eqnarray}
L_{q} = \left( \begin{array}{cc}
\frac{\displaystyle 1}{\displaystyle q} l_{q3} & l_{q-} \\
l_{q+} & - q l_{q3}
\end{array} \right)
\end{eqnarray}
and
\begin{eqnarray}
l_{q+} = \sqrt{q} k e, \qquad l_{q-} = \sqrt{q} f k, \nonumber \\
l_{q3} = \frac{1}{ [2]_{q} }( q e f - \frac{1}{q} f e )
\end{eqnarray}
It follows from (2.14), (2.17) and (1.6) 
\begin{equation}
[\hat R,(L_{q2} \hat RL_{q2}+\frac{qC_{q}}{[2]_{q}}L_{q2})]=0
\end{equation}
or in detail
\begin{eqnarray}
ql_{q3}l_{q+} - \frac{1}{q}l_{q+}l_{q3} & = & \frac{1}{[2]_{q}}C_{q}l_{q+},
\nonumber \\
\lambda l_{q3}^{2} + \frac{1}{[2]_{q}}(l_{q+}l_{q-} - l_{q-}l_{q+}) & = &
\frac{1}{[2]_{q}}C_{q}l_{q3},
\nonumber \\
ql_{q-}l_{q3} - \frac{1}{q}l_{q3}l_{q-} & = & \frac{1}{[2]_{q}}C_{q}l_{q-}
\end{eqnarray}

Having in mind the action of the $su_{q}(2)$ on the $SU_{q}(2)$ 
let us define the $SU_{q}(2)$-invariant Laplace operator by
\begin{equation}
\Delta_{q} = \frac{1}{[2]_{q}}{\rm tr}_{q}L_{q}^{2} \,. 
\end{equation}
The explicit calculation gives
\begin{equation}
\Delta_{q}=l_{q3}^{2}+\frac{1}{[2]_q}(ql_{q-}l_{q+}+\frac{1}{q}l_{q+}l_{q-})
\end{equation}
or, using (2.1)
\begin{equation}
\Delta_{q}=\frac{1}{\lambda^{2}[2]_{q}^{2}}(C_{q}+[2]_{q})(C_{q}-[2]_{q}) \,. 
\end{equation}
Due to the characteristic equation for ${\Bbb L}$ \cite{8} we also have
\begin{equation}
L_{q}^{2} + \frac{C_{q}}{[2]_{q}} L_{q} = I_{2}\otimes\Delta_{q} \,. 
\end{equation}
As in the usual $q=1$ case irreducible representations of ${\rm su}_{q}(2)$
are parametrized by the spin $l=0,\,\frac{1}{2},\,1,...$ 
The explicit formulas are \cite{9}
\begin{eqnarray}
k\mid l,m\rangle = q^{-m}\mid l,m\rangle, \nonumber \\
e\mid l,m\rangle = \sqrt{[l-m]_{q}[l+m+1]_{q}}\mid l,m+1\rangle, \nonumber \\
f\mid l,m\rangle = \sqrt{[l-m+1]_{q}[l+m]_{q}}\mid l,m-1\rangle
\end{eqnarray}
then one gets for $l_{q\pm}$ and $l_{q3}$
\begin{eqnarray}
l_{q+}\mid l,m\rangle & = & \sqrt{q^{-(2m+1)}[l-m]_{q}[l+m+1]_{q}}\mid l,m+1
\rangle,\nonumber \\
l_{q3}\mid l,m\rangle&=&\frac {1}{[2]_{q}}q^{-m}(q^{l+1}[l+m]_{q}-q^{-(l+1)}
[l-m]_{q})\mid l,m\rangle, \nonumber \\
l_{q-}\mid l,m\rangle&=& \sqrt{q^{-(2m-1)}[l-m+1]_{q}[l+m]_{q}}
\mid l,m-1\rangle
\end{eqnarray}
as well as for $C_{q}$ and $\Delta_{q}$
\begin{eqnarray}
C_{q}\mid l,m\rangle = (q^{2l+1} + q^{-(2l+1)})\mid l,m\rangle, \nonumber \\
\Delta_{q}\mid l,m\rangle = [l]_{q^{2}}[l+1]_{q^{2}}\mid l,m\rangle
\end{eqnarray}
where each multiplet $\mid l,m\rangle,\, m = -l, -l+1,...,\, l$ forms the
basis of the corresponding irreducible representation space $V_{l}$.

\newsection{$SU_{q}(2)$ and $S^{2}_{q\mu}$ as noncommutative manifolds}
The right $SU(2)_{q}$-comodule structures on $F_{q}(G)$ and $F_{q}(S)$ (1.7),
(1.17) define the corresponding left ${\rm su}_{q}(2)$-module structures on
them \cite{2},\cite{16}. It means that for every $\psi \in {\rm su}_{q}(2)$
\begin{eqnarray}
\hat \psi^{G}(u)  & = & (id\otimes \psi)\triangle (u), \nonumber \\
 & & u \in F_{q}(G)   \\
\hat \psi^{S}(v) & = & (id\otimes \psi)\varphi_{R} (v), \nonumber \\
 & & v\in F_{q}(S)
\end{eqnarray}
and the maps $\hat \psi^{G}$ and $\hat \psi^{S}$ act as linear operators
on elements of $F_{q}(G)$ and $F_{q}(S)$.
Since the mappings $\psi \rightarrow \hat \psi ^{G,S} $ are homomorphisms
\cite{5},\cite{6} all the relations (2.1), (2.2), (2.14) are valid for
$\hat k^{G,S}$, $\hat k^{-1 G,S}$, $\hat e^{G,S}$, $\hat f^{G,S}$,
$\hat L^{\pm G,S}$ and $\hat L^{G,S}$.

The elements of $F_{q}(G)$ and $F_{q}(S)$ can be also considered as left
multiplication operators on these spaces. 
Then as it was shown in \cite{5}, \cite{6} the action on the quantum group 
of the dual quantum algebra can be written as 
\begin{equation}
q\hat{\Bbb L}^{G}_{1}T_{2} = T_{2}\hat R \hat{\Bbb L}^{G}_{2}\hat R \,.
\end{equation}
Using the same technique as in one can derive the 
corresponding relation for the quantum sphere with 
$\hat{\Bbb L}^{S}$ and $M$ \cite{7}, \cite{8}
\begin{equation}
\hat R\hat{\Bbb L}^{S}_{2}\hat RM_{2} = 
M_{2}\hat R\hat{\Bbb L}^{S}_{2}\hat R \,.
\end{equation}

It is easy to see from (1.12), (2.23) and (3.4) that operator
\begin{equation}
{\cal K}=\lambda^2{\rm tr}_{q}M{\Bbb L}+\mu(1_{S}-\frac{\hat C_{q}}{[2]_{q}})
\end{equation}
commute both with $M$ and ${\Bbb L}$ and ${\cal K}(1_{S})=0$. So it is
a zero operator.

The joint algebra with the entries of $\hat{\Bbb L}^{G}$ and $T$ 
as the generators 
and relations (1.2), (2.14), (3.3) defines the algebra of functions on
quantum cotangent bundle of $SU_{q}(2)$ \cite{5}, \cite{6}. We shall denote
it by $F_{q}(T^{*}G)$. By the same arguments the joint
algebra with generators as the entries of $\hat{\Bbb L}^{S}$ and $M$ , 
relations (1.14), (2.14),(3.4)  and the additional relation:
\begin{equation}
{\cal K}=0
\end{equation} 
is called 
the algebra of functions on the quantum cotangent bundle of $S^{2}_{q\mu}$
and is denoted by $F_{q}(T^{*}S)$.

The right coaction
\begin{equation}
\varphi_{R}(T) = \triangle T, \qquad
\varphi_{R}(\hat{\Bbb L}^{G,S}) = S(T)\hat{\Bbb L}^{G,S}T, \qquad
\varphi_{R}(M) = S(T)MT
\end{equation}
define on $F_{q}(T^{*}G)$ and $F_{q}(T^{*}S)$ structures of right
$SU_{q}(2)$-comodule algebras.

Scalar products $\langle\cdot ,\cdot \rangle_{G}$ on
$F_{q}(G)$ and $\langle\cdot ,\cdot \rangle_{S}$ on $F_{q}(S)$ defines
$^{*}$-conjugation on these algebras. On $T$ and $M$ it coincides
with  the initial definition
\begin{equation}
T^{\dagger} = S(T) \,, \qquad M^{\dagger} = M
\end{equation}
and as it will be proved in Appendix the corresponding formulas for
$\hat L^{\pm}$ are
\begin{equation}
(\hat L^{\pm})^{\dagger} = S (\hat L^{\mp})
\end{equation}
(where indices $"G,S"$ are omitted). From (3.7) we have
\begin{equation}
\hat k^* = \hat k, \qquad \hat e^* = \hat f
\end{equation}
From (2.13), (2.17) and (2.19) we also have
\begin{equation}
\hat{\Bbb L}^{\dagger}=\hat{\Bbb L},\qquad \hat L_{q}^{\dagger}=\hat L_{q},
\qquad \hat C_{q}^{*} = \hat C_{q}.
\end{equation}

Now describe the $q\rightarrow 1$ limit. Let
\begin{displaymath}
\hat L^{G,S} = \left(
\begin{array}{cc}
\hat l_{3}^{G,S} & \hat l_{-}^{G,S} \\
\hat l_{+}^{G,S} & -\hat l_{3}^{G,S}
\end{array}\right)
\end{displaymath}
and $\hat C^{G,S}$ defines the $q\rightarrow 1$ limits of
$\hat L^{G,S}_{q}$ and $\hat C^{G,S}_{q}$. Then the direct calculation
using (3.1), (3.2), (2.4), (2.8b) and (2.9) gives
\begin{eqnarray}
\hat C^{G,S} & = & 2  \\
\hat l_{+}^{G} & = & a\frac{\partial}{\partial b} - b^{*}\frac{\partial}
{\partial a^{*}}, \nonumber \\
\hat l_{3}^{G} & = & \frac{1}{2}(a\frac{\partial}{\partial a} + b^{*}\frac
{\partial}{\partial b^{*}} - b\frac{\partial}{\partial b} - a^{*}\frac
{\partial}{\partial a^{*}}), \nonumber \\
\hat l^{G}_{-} & = & b\frac{\partial}{\partial a} - a^{*}\frac{\partial}
{\partial b^{*}}
\end{eqnarray}
and
\begin{equation}
\hat l^{S}_{k} = \frac{1}{i}\varepsilon_{kmn}x_{m}\frac{\partial}
{\partial x_{n}}
\end{equation}
where $\varepsilon_{kmn}$ is the Levi-Civita tensor with 
$\varepsilon_{123} = 1$ and
\begin{displaymath}
\hat l^{G,S}_{\pm} = \hat l_{1}^{G,S} \pm i\hat l^{G,S}_{2}
\end{displaymath}
Operators $\hat l^{G,S}_{k}$ satisfy the ${\rm su}(2)$ Lie
algebra commutation relations
\begin{equation}
[\hat l^{G,S}_{k},\hat l^{G,S}_{m}] = i\varepsilon_{kmn}\hat l^{G,S}_{n}
\end{equation}
and 
the $q\rightarrow 1$ limit of
${\hat \Delta}^{G,S}_{q}$ coincides with the corresponding Casimir operator
\begin{equation}
\hat \Delta = \hat l_{1}^{2} + \hat l^{2}_{2} + \hat l_{3}^{2}
\end{equation}
(where we again omitted indices $"G,S"$), and the equation (3.6) gives
\begin{equation}
\sum_{k=1}^{3}x_{k} \hat l_{k}=\frac{\displaystyle 1}{\displaystyle i}
\sum_{i=1}^{3}\varepsilon_{kmn}x_{k}x_{m}\frac{\partial}{\partial x_{n}}=0
\end{equation}

\newsection{Dirac operators on $SU_{q}(2)$ and $S^{2}_{q\mu}$}

Since all the formulas of this paragraph are similar in both 
cases of $SU_{q}(2)$ and $S^{2}_{q\mu}$ 
we shall omit indices "$G$" and "$S$".

First let us notice that according to (2.27) we can put the
$\hat L_{q}$ from (2.20) as the definition of Dirac operator on
$F_{q}(G)$ (or $F_{q}(S)$). We shall denote it by $D_{q}$.
\begin{equation}
D_{q}=\left(\begin{array}{cc}
\frac{\displaystyle 1}{\displaystyle q}l_{q3} & l_{q-} \\
l_{q+} & -ql_{q3}
\end{array}\right)
\end{equation}
The characteristic equation (2.27) gives
\begin{equation}
D_{q}^{2} + \frac{\hat C_{q}}{[2]_{q}}D_{q} = I_{2} \otimes \Delta_{q}
\end{equation}
Due to (3.11) $D_{q}^{*} = D_{q}$ and
in the case $q=1$ the corresponding operator $D$ has the form \cite{4}
\begin{equation}
D = \sum_{i=1}^{3}\sigma_{k}\otimes \hat l_{k}
\end{equation}
where $\sigma_{k}$ are the Pauli matrices
\begin{displaymath}
\sigma_{1} = \left(\begin{array}{cc}
0&1 \\
1&0
\end{array}\right), \qquad \sigma_{2} = \left(\begin{array}{cc}
0&-i \\
i&0
\end{array}\right), \qquad \sigma_{3} = \left(\begin{array}{cc}
1&0 \\
0&-1
\end{array}\right)
\end{displaymath}
and operators $\hat l_{k}$ are given by (3.12), (3.13).

The eq. (4.1) gives
\begin{equation}
D^{2} + D = I_{2}\otimes \Delta
\end{equation}
and according to \cite{10} $D$ can be interpreted as Dirac operator on
$SU(2)$ or $S^{2}$.

In order to define corresponding to $D_{q}$ spinor states we shall study
in detail the structure of the operator $D_{q}$.

Let us consider two representations of ${\rm su}_{q}(2)$: the fundamental 
irrep $\pi_{2}$ and the regular one $\pi_{reg}$ 

\begin{eqnarray}
\pi_{2}(k) = \left( \begin{array}{cc}
\frac{\displaystyle 1}{\displaystyle \sqrt{q}} & 0 \\
0 & \sqrt{q}
\end{array} \right), & &
\pi_{2}(k^{-1}) = \left( \begin{array}{cc}
\sqrt{q} & 0 \\
0 & \frac{\displaystyle 1}{\displaystyle \sqrt{q}}
\end{array} \right),\nonumber \\
\pi_{2}(e) = \left( \begin{array}{cc}
0 & 1 \\
0 & 0
\end{array} \right), & &
\pi_{2}(f) = \left( \begin{array}{cc}
0 & 0 \\
1 & 0
\end{array} \right)
\end{eqnarray}
and
\begin{equation}
\pi_{reg}(u) = \hat u
\end{equation}
for every $u\in{\rm su}_{q}(2)$. We need also their tensor product $\pi_{tot}
=\pi_{2}\otimes \pi_{reg}$ defined by
\begin{equation}
\pi_{tot}(u) = (\pi_{2}\otimes \pi_{reg})(\triangle u)
\end{equation}
for $u \in {\rm su}_{q}(2)$. Let $K = \pi_{tot}(k)$, $K^{-1} = \pi_{tot}
(k^{-1})$,
$E = \pi_{tot}(e)$ and $F = \pi_{tot}(f)$, then 
the straightforward calculation gives
\begin{eqnarray}
K & = & \left( \begin{array}{cc}
\frac{\displaystyle 1}{\displaystyle \sqrt{q}} & 0 \\
0 & \sqrt{q}
\end{array} \right) \otimes \hat k, \qquad
K^{-1} = \left( \begin{array}{cc}
\sqrt{q} & 0 \\
0 & \frac{\displaystyle 1}{\displaystyle \sqrt{q}}
\end{array} \right) \otimes \hat k^{-1},\nonumber  \\
E & = & \left( \begin{array}{cc}
0 & 1 \\
0 & 0
\end{array} \right) \otimes \hat k + \left( \begin{array}{cc}
\sqrt{q} & 0 \\
0 & \frac{\displaystyle 1}{\displaystyle \sqrt{q}}
\end{array} \right) \otimes \hat e \,,  \nonumber \\
F & = & \left( \begin{array}{cc}
0 & 0 \\
1 & 0
\end{array} \right) \otimes \hat k + \left( \begin{array}{cc}
\sqrt{q} & 0 \\
0 & \frac{\displaystyle 1}{\displaystyle \sqrt{q}}
\end{array} \right) \otimes \hat f \,.
\end{eqnarray}
Operators $K$, $K^{-1}$, $E$ and $F$ satisfies all relations (2.1) and
the corresponding central element
\begin{equation}
C^{tot}_{q} = \frac{1}{q}K^{2} + qK^{-2} + \lambda^{2}FE
\end{equation}
can be expressed in terms of $\hat C_{q}$ and $D_{q}$
\begin{equation}
C^{tot}_{q} = \frac{[2]_{q^{2}}}{[2]_{q}}I_{2}\otimes \hat C_{q}
+ \lambda^{2}D_{q}.
\end{equation}
Thus we see that $D_{q}$ due to centrality of $\hat C_{q}$ and $C^{tot}_{q}$
also commutes with $K$, $K^{-1}$, $E$ and $F$, demonstrating invariance 
property. So we may consider the
algebra ${\cal A}_{q}$ with generators $K$, $K^{-1}$, $E$, $F$ and $D_{q}$.
According to (4.1) and (2.26)
\begin{equation}
D^{2}_{q} + \frac{\hat C_{q}}{[2]_{q}}D_{q} = \frac{1}{\lambda^{2}[2]_{q}^{2}}
I_{2}\otimes (\hat C_{q} + [2]_{q})(\hat C_{q} - [2]_{q})
\end{equation}
where the operator $I_{2}\otimes \hat C_{q}$ also lies in ${\cal A}_{q}$
according to (4.8) and (4.9). Introducing the operators 
\begin{eqnarray}
{\cal L}_{+} & = & \pi_{tot} (l_{+}) = \sqrt{q} KE, \nonumber \\
{\cal L}_{3} & = & \pi_{tot} (l_{3}) = \frac{1}{[2]_{q}} (qEF -
\frac{1}{ {\displaystyle q} } FE), \nonumber \\
{\cal L}_{-} & = & \pi_{tot} (l_{-}) = \sqrt{q} FK.
\end{eqnarray}
and using (4.8), one gets 
\begin{eqnarray}
{\cal L}_{+} & = & \frac{1}{[2]_{q}}\left( \begin{array}{cc}
0 & 1 \\
0 & 0
\end{array}\right) \otimes(\hat C_{q} - \lambda [2]_{q}\hat l_{3})
+ I_{2}\otimes \hat l_{q+}, \nonumber \\
{\cal L}_{3} & = & \frac{1}{[2]_{q}^{2}} \left( \begin{array}{cc}
q & 0 \\
0 & -\frac{\displaystyle 1}{\displaystyle q}
\end{array} \right) \otimes \hat C_{q} + \frac{2}{[2]_{q}}I_{2}\otimes \hat l_{3} +
\frac{\lambda}{[2]_{q}}\left( \begin{array}{cc}
0 & \hat l_{-} \\
\hat l_{+} & 0
\end{array} \right), \nonumber \\
{\cal L}_{-} & = & \frac{1}{[2]_{q}} \left( \begin{array}{cc}
0 & 0 \\
1 & 0
\end{array} \right) \otimes (\hat C_{q} - \lambda [2]_{q}\hat l_{3})
+ I_{2}\otimes \hat l_{-}
\end{eqnarray}

We may consider now operators ${\cal L}_{+}$, ${\cal L}_{3}$, ${\cal L}_{-}$,
$\hat C_{q}$ and $D_{q}$ as generators of another $SU_{q}(2)$-covariant
algebra ${\cal B}_{q}$. Introducing the matrix generator
\begin{equation}
{\Bbb L}^{tot} = \frac{1}{[2]_{q}}I_{2}\otimes C^{tot}_{q}
+ \frac{\lambda}{\displaystyle q}L^{tot}
\end{equation}
where
\begin{equation}
L^{tot} = \left(\begin{array}{cc}
\frac{\displaystyle 1}{\displaystyle q}{\cal L}_{3} & {\cal L}_{-} \\
{\cal L}_{+} & -q{\cal L}_{3}
\end{array}\right)
\end{equation}
we can write commutation relations in the ${\cal B}_{q}$ as
\begin{equation}
\hat R^{tot}{\Bbb L}^{tot}_{2}\hat R^{tot}{\Bbb L}^{tot}_{2} =
{\Bbb L}^{tot}_{2}\hat R^{tot}{\Bbb L}^{tot}_{2}\hat R^{tot}
\end{equation}
where $\hat R^{tot} = \hat R \otimes I_{2}$.
The right $SU_{q}(2)$-coaction on ${\cal B}_{q}$ is given by (2.15) and
\begin{equation}
\varphi_{R}(D_{q}) = D_{q}\otimes I, \qquad
\varphi_{R}({\Bbb L}^{tot}) = S(T){\Bbb L}^{tot}T.
\end{equation}

The spectrum of $D_{q}$ can be easily obtained from (4.1) and (2.30). The
straightforward calculation gives two series of eigenvalues 
\begin{equation}
\lambda_{l}^{+} = [l]_{q^{2}} \,, \qquad \lambda^{-}_{l} = -[l+1]_{q^{2}} \,. 
\end{equation}
Corresponding to $\lambda^{\pm}_{l}$ eigenfunctions can be obtained
by decomposition of the tensor product $V_{\frac{1}{2}}\otimes V_{l}$ into
the direct sum $V_{l + \frac{1}{2}}\oplus V_{l - \frac{1}{2}}$
\begin{equation}
\mid l\pm \frac{1}{2},m\rangle = \sum_{m_{1},m_{2}}\left[\begin{array}{ccc}
\frac{1}{2} & l & l\pm \frac{1}{2} \\
m_{1} & m_{2} & m
\end{array}\right]_{q}
\mid \frac{1}{2},m_{1}\rangle\otimes \mid l,m_{2}\rangle
\end{equation}
where $\left[\begin{array}{ccc}
\frac{1}{2} & l & l\pm \frac{1}{2} \\
m_{1} & m_{2} & m
\end{array}\right]_{q}$ are the quantum Clebsch-Gordan coefficients \cite{12}
Using decomposition of $F_{q}(G)$ and $F_{q}(S)$ on irreducible subspaces
under the ${\rm su}_{q}(2)$ action \cite {15},\cite {16} we can express
vectors $\mid l,m_{2}\rangle$ of (4.19) in terms of $q$-special functions.

According to (2.30) spaces $V_{l\pm \frac{1}{2}}$ are eigenspaces of
$C_{q}^{tot}$ with eigenvalues $q^{2(l\pm \frac{1}{2}) + 1} +
q^{-(2(l\pm \frac{1}{2}) + 1)}$. Since in our case $V_{l\pm \frac{1}{2}}$
are imbedded into $V_{\frac{1}{2}}\otimes V_{l}$ they
also are eigenspaces of $I_{2}\otimes \hat C_{q}$
with eigenvalue $q^{2l + 1} + q^{-(2l + 1)}$. So the eq. (4.10) gives
\begin{equation}
D_{q}\mid l \pm \frac{1}{2},m\rangle = \lambda^{\pm}_{l}\mid l \pm \frac{1}{2},m\rangle.
\end{equation}

{\bf Acknowledgments.} 
One of the authors is thankful to Prof. J.A. de Azcarraga for valuable 
discussions and Generalitat Valenciana for financial support. This work 
was supported in part by RFFI Grant 96-01-00311. 

\setcounter{section}{0}\renewcommand{\thesection}{\Alph{section}}

\newappendix{}
Let us prove the general relation 
\begin{equation}
^{*}\circ \hat L^{\pm}\circ^{*} = (\hat L^{\mp})^{t}
\end{equation}
where "$^{*}$" is the $^{*}$-conjugation in $F_{q}(G)$ or $F_{q}(S)$ and "$t$"
means the usual $2\times 2$-matrix transposition.

In the $F_{q}(G)$ case (A.1) means the following
\begin{equation}
(\hat L_{1}^{\pm G}(T_{2}^{c}))^{c} = (\hat L_{1}^{\mp G})^{t}(T_{2})
\end{equation}
where $T^{c} = ST^{t}$. So from (2.10) we have $\hat L_{1}^{\pm G}(T_{2}^{c})
= S \hat L^{\pm G}_{1}(T_{2}^{t})$, and direct calculation using the 
technique of \cite{6} gives
\begin{equation}
\hat L_{1}^{\pm G}(T_{2}^{c}) = T_{2}^{c}((R^{\pm})^{-1})^{t_{2}}
\end{equation}
where $"t_{2}"$ is the transposition in the second space. For the right
hand side of (A.2) we have
\begin{equation}
(\hat L_{1}^{\pm G})^{t}(T_{2}) = T_{2}(R^{\mp})^{t_{1}}
\end{equation}
where $"t_{1}"$ means transposition in the first space. Eqs. (A.2), (A.3)
and (A.4) give
\begin{equation}
(R^{\mp})^{t_{1}} = ((R^{\pm})^{-1})^{t_{2}}
\end{equation}
or
\begin{equation}
(R^{\pm})^{-1} = (R^{\mp})^{t}
\end{equation}
and from (2.11) we get a very simple relation
\begin{equation}
\hat R^{t} = \hat R
\end{equation}
which follows immediately from the explicit expression (1.4).

In the $F_{q}(S)$ case eq. (A.1) means
\begin{equation}
(\hat L^{\pm S}(M_{2}^{t}))^{c} = (\hat L^{\mp S})^{t}(M_{2})
\end{equation}
Following \cite{6} we can easily obtain
\begin{equation}
\hat L_{1}^{\pm S}(M_{2}) = (R^{\pm})^{-1}M_{2}R^{\pm}
\end{equation}
So from (A.8) and (A.9) it follows
\begin{equation}
(R^{\pm})^{t_{2}}M_{2}((R^{\pm})^{-1})^{t_{2}} = ((R^{\mp})^{-1})^{t_{1}}
M_{2}(R^{\mp})^{t_{1}}
\end{equation}
which again leads to (A.6).

In terms of the generators eq. (A.1) means
\begin{equation}
^{*}\circ \hat k \circ ^{*} = \hat k^{-1}, \qquad
^{*}\circ \hat e \circ ^{*} = -\frac{1}{q}\hat f
\end{equation}

Let us prove now the formula (3.10). 
Consider the following chains of relations using (A.1), (2.1), 
(2.8a), (3.1) and (3.2) 
\begin{displaymath}
\langle x,\hat k(y)\rangle = h(x^{*}\hat k(y)) = h\circ \hat k(\hat k^{-1}
(x^{*})y) =
\end{displaymath}
\begin{equation}
= h\circ (id\otimes k)\triangle (\hat k^{-1}(x^{*})y) = k(h\otimes id)
\triangle (\hat k^{-1}(x^{*})y) =
\end{equation}
\begin{displaymath}
= k(I)h(\hat k^{-1}(x^{*})y) = \langle(\hat k^{-1}(x^{*}))^{*}y\rangle =
\langle\hat k(x),y\rangle
\end{displaymath}
and
\begin{displaymath}
\langle x,\hat e(y)\rangle=h(x^{*}\hat e(y))=h\circ\hat e(\hat k(x^{*})y)-
h(\hat e\hat k(x^{*})\hat k(y)) =
\end{displaymath}
\begin{displaymath}
= e\circ (h\otimes id)\triangle (\hat k(x^{*})y) - qh\circ \hat k(\hat e
(x^{*})y) =
\end{displaymath}
\begin{equation}
= e(I)h(\hat k(x^{*})y) - qh\circ (id\otimes k)\triangle(\hat e(x^{*})y) =
\end{equation}
\begin{displaymath}
= -qk(h\otimes id)\triangle(\hat e(x^{*})y) = -qh(\hat e(x^{*})y) =
\end{displaymath}
\begin{displaymath}
-q\langle(\hat e(x^{*}))^{*},y\rangle = \langle\hat f(x),y\rangle.
\end{displaymath}


\begin{thebibliography}{16}
\bibitem{1} A. Connes, ``Non-Commutative Geometry'', IHES/H/93/54.
\bibitem{2} L.D. Faddeev, N.Yu. Reshetikhin, L.A. Takhtajan,
{\it Quantization of Lie groups and Lie algebras}, Alg. 
Analiz {\bf 1} (1989) 178 (in Russian); Leningrad Math. J. {\bf 1} (1990) 193.
\bibitem{3} P. Podles, {\it Quantum spheres}, 
Lett. Math. Phys. {\bf 14} (1987) 193.
\bibitem{4}  K. Ohta and H. Suzuki, {\it Dirac operators on quantum two 
spheres}, Mod. Phys. Lett. A, {\bf 9} (1994) 2325.
\bibitem{5} A.YU. Alekseev, L.D. Faddeev, 
{\it $(T^{*}G)_{t}$: A toy model for conformal field theory}, 
Comm. Math. Phys. {\bf 141} (1991) 413. 
\bibitem{6} P. Shupp, P. Watts and  B. Zumino, {\it Bicovariant quantum 
algebras and quantum Lie algebras}, Comm. Math. Phys. {\bf 157} (1993) 305.
\bibitem{7} P.P. Kulish and R. Sasaki, {\it Covariance properties of 
reflection equation algebras}, Prog. Theor. Phys. {\bf 89} (1993) 741.
\bibitem{8} P.P. Kulish, {\it Quantum groups, $q$-oscillators and covariant
algebras}, Teor. Mat. Fiz. ( in Russian) {\bf 94} (1993) 193.
\bibitem{9} M. Jimbo, {\it A $q$-difference analogue of $U(g)$ and the
Yang-Baxter equation}, Lett. Math. Phys. {\bf 10} (1985) 63.
\bibitem{10} N.Berline, E. Getzler, M. Vergne, "Heat Kernels and Dirac
Operators." Berlin, 1992.
\bibitem{11} L.L.Vaksman and Ya.S. Soibelman, {\it Algebra of functions
on quantum group $SU(2)$}, Funct. Anal. Pril. {\bf 22} (1988) 1.
\bibitem{12} A.N. Kirillov and N.Yu. Reshetikhin, {\it Representations
of the algebra $U_{q}(sl(2))$, $q$-orthogonal polynomials and invariants
of links}, preprint LOMI E-9-88 (1988); Infinite-dimensional Lie algebras and
groups, W. S., Singapore (1989).
\bibitem{13} J.A. de Azcarraga, P.P. Kulish, F. Rodenas, {\it Quantum
groups and deformed special relativity }, Fortschr. der Phys. {\bf 44} (1996) 
1; hep-th/9405161.
\bibitem{14} H. Grosse, C. Klimcik, P. Presnajder, {\it Simple field
theoretical models on noncommutative manifolds}, preprint CERN-TH/95-138; 
hep-th/9510177; hep-th/9510083; hep-th/9505175.
\bibitem{15} M.Noumi and K. Mimachi, {\it Rogers's $q$-ultraspherical
polynomials on a quantum $2$-sphere}, Duke Math. Jour. {\bf 63} (1991) 65.
\bibitem{16} T. Masuda, K. Mimachi, Y. Nakagami, M. Noumi, K. Ueno,
{\it Representations of the quantum group $SU_{q}(2)$ and the little $q$-Jacobi
polynomials}, Journ. Func. Anal. {\bf 99} (1991) 357.
\end{thebibliography}
\end{document}